\documentclass[10pt]{article}
\usepackage{amsmath, amssymb}
\usepackage{amscd}

\newtheorem{theor}{Theorem}[section]

\newtheorem{lem}[theor]{Lemma}
\newtheorem{propo}[theor]{Proposition}
\newtheorem{coro}[theor]{Corollary}

\newenvironment{pf}{{\it Proof.}}{\hfill $\square$\\}

\begin{document}
\title{\bf Dynamics on the space
of harmonic functions and  
the foliated Liouville problem}
\author{R. Feres and A. Zeghib}
\maketitle
      
\footnotetext{Dept. of Mathematics - 1146, Washington University,
 St. Louis, MO 63130, USA.}
 \footnotetext{UMPA - \'Ecole Normale Sup\'erieure de Lyon,
69364 Lyon CEDEX 07, France.}      

\begin{abstract}
We study here the  action of subgroups of $PSL(2,\mathbb{R})$
on the space of harmonic functions on the unit disc bounded by a common constant,
as well as the relationship this action has
  with the foliated Liouville problem:
Given a foliation of a compact manifold by Riemannian leaves
and a leafwise harmonic continuous function on the manifold, is the function
leafwise constant? We give  a number of positive results
and also show a general class of  examples for which the Liouville property
does not hold. 
 The connection between the Liouville property and the dynamics on
 the space of harmonic functions as well as general properties
of this dynamical system are explored. It is shown among other properties that
the $\mathbb{Z}$-action generated by hyperbolic or parabolic elements of 
$PSL(2,\mathbb{R}) $ is chaotic.
\end{abstract}

\section{Introduction}
Let $\text{Har}(\mathbb{D})$ denote the space of complex valued harmonic functions on the unit disc
   $\mathbb{D}=\{z\in \mathbb{C}:|z|<1\}$ bounded by a common constant,
which will be taken without loss of generality to be $1$.
Endowed with the topology of uniform convergence on
compact subsets of $\mathbb{D}$, $\text{Har}(\mathbb{D})$ is a compact metrizable space.
 The group $G=PSL(2,\mathbb{R})$ of  hyperbolic 
 isometries of $\mathbb{D}$ acts on $\text{Har}(\mathbb{D})$
by composition: $(g, f)\mapsto f\circ g^{-1}$, where $f\in  \text{Har}(\mathbb{D})$ and
$g\in G$.
A central concern of the present paper is   the dynamics of the action
 on $\text{Har}(\mathbb{D})$ of  $G$ and its subgroups, in particular, of
lattices in $G$.

We are led to study the dynamics of lattices in $G$ on $\text{Har}(\mathbb{D})$ by
what will be call  here the {\em foliated Liouville problem}.
The general setting for this problem is 
a foliated space $(M,\mathcal{F})$ (as defined in \cite{candelconlon}),
 where $M$ is a compact topological
space and the leaves of   $\mathcal{F}$ are  smooth Riemannian manifolds.
The Riemannian metric  and  all its derivatives  are
  assumed to vary continuously on $M$.

The foliated space
 $(M,\mathcal{F})$ will be said to have the  {\em Liouville
property}
if   continuous leafwise harmonic  functions on $M$ are leafwise constant.
More precisely, let $\Delta$ denote the tangential Laplace operator
associated to the Riemannian metric on leaves of $\mathcal{F}$. 
A continuous 
function $f$ which is  smooth along leaves   and satisfies $\Delta f=0$
will be called {\em leafwise harmonic}. Then the Liouville
property holds if for any such   function its
restriction to each leaf of $\mathcal{F}$ is constant. 
 
If the leaves of $\mathcal{F}$  are
complex rather than Riemannian manifolds, 
the related problem of deciding when leafwise {\em holomorphic}    continuous functions
on $M$ are leafwise constant was studied
in our paper \cite{ghani}. 
Clearly, these two problems are related. For example, if the leaves are
K\"ahler manifolds (as in the main class of examples of
foliated bundles over Riemann surfaces), the real and imaginary parts of 
holomorphic functions are harmonic, so that some of the results 
obtained in our earlier paper have immediate implications to the present setting.
 On the other hand, some tools used in that paper are
unique to holomorphic functions (such as the open mapping principle),
which makes the harmonic case more difficult and more interesting.
 Since both the harmonic and holomorphic versions of the problem will arise
  in the course of the present paper, it will be convenient at times to use the alternative terminology: 
a foliated space with the Liouville property will also be called {\em harmonically simple},
while  its holomorphic counterpart will be called  {\em holomorphically simple}.
(The reader should note that the latter was called in \cite{ghani} {\em holomorphically plain}.)

 This   property  obviously holds
for the trivial foliation, consisting of
a single leaf ($M$  itself). 
Less trivial, but for our concerns equally uninteresting examples are provided by
 foliations   whose leaves are in a certain sense parabolic,
so that 
  Liouville theorems    extensively studied in 
geometry and complex analysis can be  applied to each leaf separately.

When the foliation has leaves that individually admit nonconstant bounded harmonic functions,    dynamical considerations must come into play in
order to decide whether or not the Liouville property holds.   To give a rather
simple illustration of this point, note that a minimal foliation (that is, one
that does not contain proper closed saturated sets) has the Liouville property
regardless of the geometry of the leaves, as an application of the
maximum principle for harmonic functions immediately shows. More interesting related results
will be offered later. On the other hand, simply having complicated  
transversal dynamics is not by itself an obstacle to the existence of nonconstant,
leafwise harmonic continuous functions on $M$, as will be seen later
with  an example  of a codimension $2$
ergodic foliation for which the Liouville property fails.

In this paper we describe   results and examples connecting
the Liouville property  with the foliation's  transverse dynamics. 
We will consider in some detail 
foliated bundles over compact
 Riemann surfaces
of genus $2$ or greater. This is   a  particularly well suited
class of examples 
for this study since typical
leaves    
 often admit   nonconstant bounded harmonic functions, 
and such examples can be constructed so as to possess a wide range of dynamical properties.

It   will be shown that constructing  foliations without the Liouville
property,  in the class of foliated bundles over compact Riemann surfaces,   
naturally leads to study the dynamics of uniform lattices in $G$  on  $\text{Har}(\mathbb{D})$.  Among other properties, it will
be seen that   this action is topologically
transitive, and that the $\mathbb{Z}$-action generated by any hyperbolic element
in the lattice is chaotic in the sense of \cite{devaney}.

The plan of the paper is as follows. After setting some notation, we describe a number of results, topological and measure theoretic,
proving the Liouville property for certain classes of foliations. The topological results
are harmonic counterparts to results proven in \cite{ghani} in the holomorphic setting,
and most proofs, although not all, are similar to those in our previous paper.
The  measure theoretic results are mainly aimed at connecting
  our problem with L. Garnett's theory of harmonic measures for foliations.

Next, we describe a class of examples of harmonically non-simple foliations.
These are foliated bundles over a compact Riemann surface, of codimension  $2$,
and ergodic.
  It is also shown how foliated bundles over a compact Riemann surface
for which the Liouville property does not hold (and 
nontrivial leafwise harmonic functions) can be obtained from a
universal foliated bundle construction.

 We then compare the
harmonic problem with the corresponding holomorphic problem studied in
\cite{ghani}.  It is remarked that each continuous leafwise harmonic
function defines a cohomology class in the first tangential de Rham cohomology
space of a leafwise K\"ahler  foliation, and that this class vanishes if and only if
the function is the real part of a leafwise holomorphic function. The example
of a harmonically non-simple foliation referred to above has, in particular,
non-trivial first cohomology 
as it is shown that the example is holomorphically simple.

Finally, we prove a number of results about the dynamics of   actions
of subgroups of $PSL(2,\mathbb{R})$ on $\text{Har}(\mathbb{D})$.
One of the main results is that the $\mathbb{Z}$-action
of hyperbolic or parabolic elements in $PSL(2,\mathbb{R})$
define chaotic dynamical systems in the sense of \cite{devaney}.

The first author would like to thank the \'Ecole Normale Sup\'erieure de Lyon, UMPA,
for its support and hospitality while this works was being written.
\section{Notations and general facts}
We first set  some notation.  
The  unit disc $\mathbb{D}=\{z\in \mathbb{C}:|z|<1\}$  
 can be viewed alternatively as a Riemannian manifold with the
metric of constant curvature $-1$, or as a complex manifold.
In the first case, the set  of all complex valued   harmonic functions   $f$ on 
$\mathbb{D}$
such that $f(\mathbb{D})\subset \overline{\mathbb{D}}$ will 
be written  $\text{Har}(\mathbb{D})$, while $\text{Hol}(\mathbb{D})$ will
denote the set of all holomorphic functions from  $\mathbb{D}$ into $\overline{\mathbb{D}}$.   
Notice that $\text{Har}(\mathbb{D})\subset \text{Hol}(\mathbb{D})$. On the other hand, a holomorphic map whose real
part lies in $\text{Har}(\mathbb{D})$ may fail to be bounded.

The subset of  $\text{Har}(\mathbb{D})$ (resp.,
$\text{Hol}(\mathbb{D})$)  of nonconstant functions will be denoted by
$\text{Har}(\mathbb{D})^*$ (resp.,
$\text{Hol}(\mathbb{D})^*$).

Endowed  with the topology of 
uniform convergence on compact subsets of $\mathbb{D}$, both $\text{Har}(\mathbb{D})$ and $\text{Hol}(\mathbb{D})$ are compact metrizable spaces.
Recall that
$\text{Har}(\mathbb{D})$  can   be identified, by means of
the Poisson representation formula, with the unit ball
of $L^\infty(S^1)$ endowed with the weak*-topology, where $S^1=\partial \mathbb{D}$. Thus, denoting
the Poisson kernel by
$K(z,\theta)=(1-|z|^2)/(|e^{2\pi i\theta}-z|^2),
$ 
convergence in $\text{Har}(\mathbb{D})$
\begin{equation*}\varphi_n(\cdot)=\int_0^1K(\cdot,\theta)f_n(\theta)d\theta\rightarrow
\varphi(\cdot)=\int_0^1K(\cdot,\theta)f(\theta)d\theta
\end{equation*}  holds  if  for every $g\in L^\infty(S^1)$ (equivalently, for  each    $g$ in a dense subset
of  continuous functions in $L^1(S^1)$),    one has
\begin{equation*}\int_0^1(f_n(\theta)-f(\theta))g(\theta)d\theta\rightarrow 0
\end{equation*}
 as $n\rightarrow \infty$.

It is well known that $PSL(2,\mathbb{R})$, acting
on  $\mathbb{H}=\{z=x+iy:y>0\}$ by partial linear transformations,
$gz=(az+b)/(cz+d),$ is the group of isometries of the hyperbolic plane. The conformal isomorphism 
  $\phi:\mathbb{D}\rightarrow\mathbb{H}$ given by $\phi(z)=i(1-z)/(1+z)$, 
conjugates the action on the upper half-plane to an action 
of $PSL(2,\mathbb{R})$  on the unit disc, according to  the expression
$gz= {(\alpha z + \beta)}/{( \bar{\beta}z +\bar{\alpha})}$, where
$2\alpha={a+d} +{(b-c)} i, 2\beta= {-a+d}  +  {(b+c)}i$ and 
  $z\in \mathbb{D}$. Conjugation by $\phi$ gives, in fact,
an isomorphism between $PSL(2,\mathbb{R})$ and $PSU(1,1)$. Elements of the latter
group
are represented by complex  matrices 
$\begin{pmatrix} \alpha & \beta \\ 
\bar{\beta} & \bar{\alpha}\end{pmatrix}$ of determinant $1$, modulo the center.
The group $PSL(2,\mathbb{R})$, which consists of harmonic  (resp.,
holomorphic) self-maps of $\mathbb{D}$, naturally acts on $\text{Har}(\mathbb{D})$
(resp., $\text{Hol}(\mathbb{D})$) by:  $(g, f)\mapsto f\circ g^{-1}$, where 
$g^{-1}(z)=(\bar{\alpha}z-\beta)/(\alpha -\bar{\beta}z).$

We now recall the construction of a foliated bundle over a Riemann surface.
 Let $S$ be a compact connected
Riemann surface, $\tilde{S}$ its universal covering space, and $\Gamma$ the 
fundamental group of
  $S$. The covering action of $\gamma\in \Gamma$ on $p\in \tilde{S}$   will be written $p\gamma$.
 Let $X$ be a compact connected space on which
$\Gamma$ acts by homeomorphisms, for a given homomorphism 
$\rho:\Gamma\rightarrow \text{Homeo}(X)$ from $\Gamma$ into the group of
homeomorphisms of $X$.   
Given $\gamma\in \Gamma$ and $x\in X$,  
$\rho(\gamma)(x)$ will be written simply as $\gamma(x)$.
 $\Gamma$ naturally acts on the product $\tilde{S}\times X$, properly discontinuously, by
$
(p,x)\cdot \gamma:= (p\gamma, \gamma^{-1}(x)).$
The space of orbits for this action,
$M=(\tilde{S}\times X)/\Gamma$,  is a foliated space whose
leaves  are  transverse to the fibers of the natural projection
$\pi:M\rightarrow S$. The restriction of the projection map to individual leaves
are covering maps. The resulting foliated space will be written $(M_\rho, \mathcal{F}_\rho)$.

\section{Harmonically simple foliations}
We describe in this section conditions under which a foliated space has the Liouville
property. These are mostly results that were given in \cite{ghani} for
the holomorphic Liouville problem, but which also hold in the harmonic setting.
The proofs are essentially the same as  for the corresponding 
results in that paper, except for Proposition \ref{grom} below.
On the other hand, 
not all that what was shown in our earlier paper seems to have
an easy translation to the present setting.
For example, in Theorem 1.15 of \cite{ghani}, it is shown
that $(M,\mathcal{F})$ is holomorphically simple
whenever it  has codimension $1$. We do not know whether the 
counterpart  for harmonic functions holds.

Some of the results in this section are of a measure theoretic nature while others
are purely topological. In all cases, $(M,\mathcal{F})$ is a compact foliated space  with leafwise Riemannian metric.

\subsection{Topological results}
Let $f$ be
  a continuous leafwise harmonic function  on $M$.
It is clear that the set of leaves where $f$ is constant is a compact nonempty $\mathcal{F}$-saturated set. Notice that it is nonempty
by the maximum value property of harmonic functions, since on a leaf containing a 
point where $f$ attains a maximum or a minimum value, $f$ must be constant.
 This remark immediately
yields the following proposition. 
\begin{propo}\label{minimal}
If the closure of every leaf of $\mathcal{F}$ contains a unique minimal set,
then the leafwise Riemannian compact foliated space$(M,\mathcal{F})$
is harmonically simple.
\end{propo}
\begin{pf}
Under these assumptions, the maximum and minimum values of the function 
must coincide on each leaf closure.
\end{pf}

\begin{coro}\label{coco}
Let $(M_\rho,\mathcal{F}_\rho)$ be a  foliated bundle over a compact Riemannian manifold $S$ with compact fiber a differentiable manifold $V$, where 
$\rho:\Gamma\rightarrow H$ is a homomorphism from
the fundamental group of $S$ into a compact group of diffeomorphisms of $V$.  Then 
$(M_\rho,\mathcal{F}_\rho)$ is harmonically simple.
\end{coro}
\begin{pf}
Since $H$ is compact, there exists an $H$-invariant Riemannian metric on $V$,
making $(M_\rho,\mathcal{F}_\rho)$ a Riemannian foliation. By \cite{molino},
the closure of each leaf is a minimal set.
\end{pf}

Denote by $P(W)$ the projective space associated to a finite dimensional
real or complex vector space $W$. The general linear group  $GL(W)$ naturally acts 
on $P(W)$. 
We consider next foliated bundles with fiber $P(W)$ over a compact Riemannian 
manifold $S$, 
associated to homomorphisms $\rho:\Gamma\rightarrow GL(W)$,
where $\Gamma$ is the fundamental group of $S$.
An element $\gamma\in \Gamma$ is called {\em proximal} if the maximal
characteristic exponent of $\rho(\gamma)$ is simple. Proposition \ref{minimal}
and \cite[3.4 and 3.6, ch. VI]{margulis} yields the next proposition.

\begin{propo}\label{propA}
Let $S$ be a connected, compact, Riemannian manifold  with fundamental group $\Gamma$, $W$ a finite dimensional   vector space and 
suppose that  $\rho:\Gamma\rightarrow GL(W)$ is a homomorphism such that $\Gamma$ contains a proximal element.
Let $(M_\rho,\mathcal{F}_\rho)$ be the foliated bundle over $S$ with fiber
$P(W)$ and $\Gamma$-action on $P(W)$ given by $\rho$. Then 
$(M_\rho,\mathcal{F}_\rho)$ is harmonically simple.
\end{propo}

\begin{coro}\label{caca}
Let $S$, $\Gamma$,   $W$, and $(M_\rho, \mathcal{F}_\rho)$ be as in Proposition \ref{propA},
where the representation
  $\rho:\Gamma
\rightarrow GL(W)$ is now such that  the projection of   $\rho(\Gamma)$ into $PGL(W)$
is Zariski dense. Then $(M_\rho, \mathcal{F}_\rho)$
is harmonically simple.
\end{coro}
\begin{pf}
If the image of $\rho(\Gamma)$ in $PGL(W)$ is not precompact, the hypothesis of
Proposition \ref{propA} hold by  \cite[Theorem 4.3(i)]{margulis}. Otherwise,
the conclusion follows from Corollary \ref{coco}.
\end{pf}

The proof of the next theorem is, with obvious changes, the same as
for Theorem 1.11 in \cite{ghani}. It uses Corollary \ref{caca}, as well as
\cite{rapin} and \cite[8.2]{burger}.
\begin{theor}
Let $S=\mathbb{D}/\Gamma$ be a compact Riemann surface
of genus $g\geq 2$
and write $G=GL(n,\mathbb{C})$. Then there is a Zariski dense
subset $U$ of the representation variety $\text{Hom}(\Gamma,G)$
such that for each $\rho\in U$, the foliated space $(M_\rho, \mathcal{F}_\rho)$
with fiber $P(\mathbb{C}^{n})$ is harmonically simple.
\end{theor}

Another class of foliations to which the ideas of this section apply 
are obtained from actions of Gromov-hyperbolic groups on their
boundary. (See \cite{bow} for the general definitions and facts concerning
such groups.)
The proof of the next theorem
requires some modifications of the proof of  its holomorphic counterpart,
  Proposition 1.13 of \cite{ghani}.

\begin{theor}\label{grom}
Let $\Lambda$ be  a  Gromov-hyperbolic group, $X$ the boundary of $\Lambda$, and $S$
a compact connected Riemannian manifold with fundamental
group $\Gamma$. Suppose that $\Gamma$ acts on $X$ via a homomorphism
$\rho:\Gamma\rightarrow \Lambda$ and let $(M_\rho,\mathcal{F}_\rho)$ be the corresponding
foliated bundle over $S$. Then $(M_\rho,\mathcal{F}_\rho)$ is harmonically simple.
\end{theor}
\begin{pf}
Let $L$ denote the limit set for the action of $\rho(\Lambda)$ on $X$.
There are only three possibilities for the cardinality of $L$: $1, 2, \infty$.
If the cardinality is either $1$ or $\infty$, the action 
of $\Lambda$ on $X$ has a  unique minimal set (see \cite{gromovhyp}), so
the conclusion of the theorem follows from  Proposition \ref{minimal}.

So suppose that $L$ consists of exactly $2$ points. In this case $\rho(\Gamma)$
fixes the geodesic joining them. But the stabilizer of a geodesic in 
a (countable) hyperbolic group is virtually cyclic. If $\rho(\Gamma)$ is
finite, all leaves of $\mathcal{F}$ are compact, so the conclusion of
the theorem holds. If  $\rho(\Gamma)$ is not finite, any noncompact
leaf is a finitely generated
abelian Galois covering of a common compact Riemann surface. But it
is known that such surfaces, individually,   satisfy the Liouville property. (See
\cite{babi}.) Therefore the theorem follows.  
\end{pf}

We note that, in the next theorem, a  non-discrete $\rho(\Gamma)$ is
not excluded.
\begin{theor}\label{a}
The previous theorem still holds after replacing $\Lambda$ by
$SL(2,\mathbb{C})$. 
\end{theor}
\begin{pf}
The proof is essentially the same as for Theorem \ref{grom}.
Of the three possibilities for the  
limit set of the action of $\rho(\Gamma)$ on the boundary of the symmetric space,
which is now $\mathbb{H}^3$,  the difficult case corresponds to
a two point limit set.
So it can be assumed that $\Gamma$ fixes a geodesic
in $\mathbb{H}^3$.
The stabilizer in $SL(2,\mathbb{C})$ of this geodesic  is $\mathbb{R}\times SO(2)$,
an abelian group. Consequently, non-compact leaves are (virtually)  
 finitely generated abelian Galois coverings of a same compact
Riemann surface. As before, \cite{babi} finishes the proof.
\end{pf}

It seems plausible to expect the statement of Theorem \ref{a}  still to be
valid if one replaces $SL(2,\mathbb{C})$ with  a general
rank-$1$ semisimple Lie group.
The proof given above, however, breaks down since the
stabilizer of a geodesic in the corresponding symmetric space
is $\mathbb{R}\times K$, for a compact  $K$ which may be non-abelian.
Nevertheless, we expect that a more detailed analysis
should yield the result in this more general case.

 \subsection{Measure theoretic results} The next  theorem is from \cite{garnett}. We recall that a measure $m$ on $M$
is said to be  {\em harmonic} if for all continuous leafwise smooth functions
$h:M\rightarrow \mathbb{R}$,
\begin{equation*}
\int_M(\Delta h)(x) dm(x)=0.
\end{equation*}
A harmonic measure on  $(M,\mathcal{F})$ is called {\em totally invariant}
if, together with the leafwise Riemannian measure, it yields a transverse
invariant measure for $\mathcal{F}$.
We refer to \cite{garnett} or \cite{candel} for the general properties and
results concerning harmonic measures.
\begin{theor}[Garnett]
Let $m$ be a harmonic measure on $(M,\mathcal{F})$ and $f$ a measurable,
$m$-integrable,
leafwise harmonic function on $M$. Then $f$ is constant on $m$-a.e. leaf.
\end{theor}

Some generalizations of Garnett's result are obtained by S. Adams in \cite{adams}.

\begin{coro}
If the union of the supports of   harmonic measures is all of $M$, then $(M,\mathcal{F})$
is harmonically simple.
\end{coro}

\begin{coro}
$(M,\mathcal{F})$ is harmonically simple whenever
it  admits a transverse invariant measure of full support.
\end{coro}

The last corollary shows that it is easy to obtain examples of 
harmonically simple foliations.
In fact, for every volume preserving action of  a cocompact lattice $\Gamma\in PSL(2,\mathbb{R})$ on 
a compact manifold $X$, the corresponding foliated bundle
$(\mathbb{D}\times X)/\Gamma$ has that property.

 \section{Harmonically nonsimple  foliations}
 We gave in \cite{ghani} the following  example of a foliated bundle $(M,\mathcal{F})$   that is not holomorphically
or harmonically
simple. 
Here
$M=(\mathbb{D}\times \mathcal{C})/\Gamma$, where
 \begin{equation*}
\mathcal{C}=\{[z_1,z_2,t]\in \mathbb{R}P^4: |z_1|^2-|z_2|^2=t^2\}
\end{equation*}
is   an $SU(1,1)$-invariant submanifold of projective space for an action
  of $SU(1,1)$ on $\mathbb{R}P^4$ defined as
follows:
\begin{equation*}
\begin{pmatrix}
\alpha & \beta\\
\bar{\beta} & \bar{\alpha}
\end{pmatrix}\cdot [z_1, z_2, t]=[\alpha z_1 +\beta\bar{z}_2, \alpha z_2 +\beta \bar{z}_1, t].
\end{equation*}
The function $f:\mathbb{D}\times \mathcal{C}\rightarrow \mathbb{C}$ given
by
\begin{equation*}
f(z, [u_1, u_2, t])=\frac{\bar{u}_1z -u_2}{-\bar{u}_2 z +u_1}
\end{equation*}
is easily shown to pass to the quotient by $\Gamma$, yielding  a real analytic
leafwise holomorphic function.

The transverse dynamics of the above example   ``essentially''
corresponds to the action of $\Gamma$ on a $3$-sphere obtained
by compactifying $PSU(1,1)$ (on which $\Gamma$ acts by translations)
 by adding a circle at infinity. It  is not a particularly
complicated dynamics, topologically or measure theoretically.

It is natural to ask whether    
dynamical   properties such as recurrence and ergodicity might somehow   prevent   examples of this kind.
The purpose of this section is to construct a  {\em harmonically} non-simple
foliated space which is ergodic with respect to a transverse measure class that is
positive on open sets. (It will be shown later that the example is holomorphically
simple.)
More precisely, we have the following theorem.

\begin{theor}\label{mainexample}
 Let $S=\mathbb{D}/\Gamma$ be a compact Riemann surface,
where $\Gamma$ is a cocompact lattice
 in  $PSL(2,\mathbb{R})$.  Then there exists
a  foliated bundle $(M,\mathcal{F})$ over $S$ with fiber $S^2$
such that:
\begin{enumerate}
\item $(M,\mathcal{F})$ is real analytic in the complement of a pair of compact
leaves, $S_1, S_2$ homeomorphic to $S$;
\item $(M,\mathcal{F})$ is ergodic with respect to the smooth measure class;
\item The   complement of $S_1\cup S_2$
  has a real analytic compactification, which 
  is  an ergodic  foliated bundle over $S$ with fiber $S^1\times [0,2\pi]$.
\item For both $(M,\mathcal{F})$ and its analytic compactification the Liouville property
does not hold.  Moreover, a continuous leafwise harmonic, not leafwise constant,
can be found that is real analytic in the complement of $S_1\cup S_2$. 
 \end{enumerate}
\end{theor}

We need to following definitions.
By an {\em arc} in $S^1$ we   refer to a set of the form
$\{\zeta e^{it\theta}: 0\leq t\leq 1\},
$
where $\zeta\in S^1$ and $\theta\in (0,2\pi)$. Therefore the space of arcs is  parametrized by the cylinder $S^1\times (0,2\pi)$ and has a natural completion,
$S^1\times [0,2\pi]$, in which  $S^1\times\{0\}$ consists of trivial (single point) arcs
and $S^1\times \{2\pi\}$ consists of full circle arcs with the  initial point specified.
Since $PSL(2,\mathbb{R})$ acts on the boundary of the unit disc, it also
acts on the space of arcs $S^1\times (0,2\pi)$. The action can be written explicitly
as follows: $g(\zeta, \theta)=(g(\zeta), \Theta_g(\zeta,\theta))$, where
$g(\zeta)=(\alpha \zeta +\beta)/(\bar{\beta}\zeta +\bar{\alpha})$, $|\alpha|^2-|\beta|^2=1$,
and
\begin{equation*}
\Theta_g(\zeta,\theta)=\int_0^1\frac{dt}{|\bar{\alpha} + \zeta \bar{\beta} e^{it\theta}|^2}.
\end{equation*}
This is  a real analytic action. It 
extends an action on  $S^1\times [0,2\pi]$ which is also real analytic. The boundary components $S^1\times \{0\}$
and $S^1\times \{2\pi\}$ are invariant subsets. 

\begin{lem}\label{lem1}
Let $\Gamma$ be an arbitrary  lattice in $PSL(2,\mathbb{R})$. Then the action of $\Gamma$
on $S^1\times [0,2\pi]$ defined above    is ergodic with respect to the smooth measure
class on the cylinder.
\end{lem} 
\begin{pf}
It is easily seen that   $PSL(2,\mathbb{R})$ acts transitively on the space 
of (nontrivial) arcs and that
the isotropy group of the arc specified by $\zeta=1$ and $\theta=\pi$
is the subgroup   $A\subset PSL(2,\mathbb{R})$ represented by diagonal matrices in $SL(2,\mathbb{R})$. Thus the action of $\Gamma$ on the space of
arcs is analytically conjugate to the action of $\Gamma$ by left translations
 on the quotient
$PSL(2,\mathbb{R})/A$. But it is well known that this action is  ergodic with respect
to the   smooth measure  class on the quotient. In fact,
its dual action,
  of $A$ on $PSL(2,\mathbb{R})/\Gamma$,
is the geodesic flow 
  of  a compact negatively curved
surface.
\end{pf}

Write $X_0=S^1\times [0,2\pi]$ and 
define   $F:\mathbb{D}\times X_0\rightarrow \mathbb{R}$ according to 
\begin{equation*}
F(z, \zeta, \theta)=\frac{\theta}{2\pi}\int_0^1\frac{1-|z|^2}{|\zeta e^{it\theta} -z|^2}dt.
\end{equation*}
Observe that $F$ is real analytic and $F(z, \zeta,0)=0$, $F(z,\zeta,2\pi)=1$, for
all $z\in \mathbb{D}$ and $\zeta\in S^1$.
\begin{lem}\label{lem2}
 We have $F(g(z), g(x))=F(z, x)$ for
all $g\in PSL(2,\mathbb{R})$, $z\in \mathbb{D}$, and $x\in  X_0$.
\end{lem}
\begin{pf} Denote by $\chi_I$ the indicator function of the arc $I$ determined by
$\zeta$ and $\theta$.
By the Poisson formula, 
 $z\mapsto \varphi(z)=F(z,\zeta,\theta)$  is the unique 
harmonic function on the unit disc with boundary value $\chi_I$.
Since the action of $PSL(2,\mathbb{R})$ on the unit disc is conformal,
$\varphi\circ g$ is the unique harmonic function on $\mathbb{D}$
with boundary value $\chi_I\circ g=\chi_{g^{-1}(I)}$, from which the claim follows.
\end{pf}

Denote by $(M_0, \mathcal{F}_0)$ the foliated bundle for which $M_0=(\mathbb{D}\times X_0)/\Gamma$ and $\mathcal{F}_0$ is the resulting
foliation by coverings of  $\mathbb{D}/\Gamma$.
Notice that 
$M_0$ is a compact manifold with boundary, whose   boundary components
are the  two  compact leaves of $\mathcal{F}_0$ associated to the
invariant boundary circles of $X_0$. 
Due to Lemma \ref{lem1}, 
$\mathcal{F}_0$ is   ergodic.

Due to Lemma \ref{lem2}, $F$ yields on the quotient $M_0$ a  leafwise harmonic
continuous 
function,  $F_0:M_0\rightarrow \mathbb{R}$. 
Notice that $F_0$ is constant on the boundary leaves.

Denote by $X_1$ the sphere $S^2$ obtained as the quotient of
$S^1\times [0,2\pi]$ by collapsing the two boundary circles to the north and
south poles.  By repeating the above construction, now with $X_1$,
we obtain a continuous foliated bundle over $\mathbb{D}/\Gamma$
with transversal fiber $S^2$, which is ergodic and admits a continuous leafwise
harmonic function.

Notice that the $\Gamma$-action on $S^1\times [0,2\pi]$ 
that was used to define $(M,\mathcal{F})$ passes to the quotient 
so as to define an action on the torus $S^1\times S^1$. The latter action
is real analytic and ergodic with
respect to the smooth measure class, so it induces an ergodic real analytic foliated bundle. This foliated bundle, however, is harmonically simple.

 We end the section with the following  two remarks. The first one  is
an immediate consequence of Garnett's theorem.
\begin{propo}
The support of any harmonic measure for the foliation $(M,\mathcal{F})$  of Theorem \ref{mainexample} is contained in the union of the compact leaves $S_1$ and $S_2$.
\end{propo}

\begin{propo}\label{holsimple}
  $(M,\mathcal{F})$,  of Theorem \ref{mainexample}, is holomorphically simple.
\end{propo}
\begin{pf}
The claim is a corollary of Proposition 1.5 of \cite{ghani}. We give  below
a simple
direct proof.
Let $f:M\rightarrow \mathbb{C}$ be a leafwise holomorphic function.
Clearly $f$ is constant on the compact leaves $S_1$ and $S_2$.
Let $K$ be the   union of leaves on which $f$ is constant.
The transversal dynamics is such that any leaf which is not $S_1$ or $S_2$
approaches $S_1$ or $S_2$ along the direction of a hyperbolic element of $\Gamma$.
Therefore $f$ can
take at most two values on $K$.
On the other hand, the image of the complement of $K$ under $f$
is an open subset of $\mathbb{C}$, due to the open mapping principle
for holomorphic functions. Since the union of this open set with
the (at most two points of the) image of $K$ is compact, it follows that
the complement of $K$ is empty. 
\end{pf}

\section{A universal non-Liouville foliation}
Foliated bundles over a compact Riemann surface of genus at least $2$
for which the Liouville property fails afford a general
description, which is explained in this section. A similar construction for
the holomorphic case was given in \cite{ghani}.

Write $X_0=\text{Har}(\mathbb{D})$, and
let  $S=\mathbb{D}/\Gamma$ be a compact Riemann surface.
As already noted, $X_0$ is a compact metrizable
space with the topology of uniform convergence on compact subsets of the unit disc,
and supports a continuous action of $PSL(2,\mathbb{R})$ given by
$(g,f)\mapsto f\circ g^{-1}$. 
Thus, it makes sense to form the compact foliated space
$M_0=(\mathbb{D}\times X_0)/\Gamma$ over $S$. We denote this 
foliated bundle by $(M_0,\mathcal{F}_0)$. 

A leafwise  harmonic
function on $M_0$ is given by the following tautological construction.
Define $\bar{\phi}:\mathbb{D}\times X_0\rightarrow  \mathbb{C}$ by $\bar{\phi}(z, f):=
f(z)$. It is easily checked that $\phi$ passes to the quotient and
defines a continuous, leafwise harmonic function on $M_0$.

Given foliated bundles 
$(M',\mathcal{F}')$ and $(M, \mathcal{F})$ over $S$, 
we define a {\em harmonic morphism}   between them as
   a continuous, fiber preserving map
$F: M'\rightarrow M$ that  sends leaves of $\mathcal{F}'$ to
 leaves of $\mathcal{F}$ and the
restriction to each leaf of $\mathcal{F}'$ is a harmonic map.

In trying to construct examples of harmonically nonsimple foliations, it
is useful to have in mind the following easy fact.
\begin{propo}
Let $(M,\mathcal{F})$ be a foliated bundle over $S=\mathbb{D}/\Gamma$ with
fiber $V$. Then there is a one-to-one  correspondence between continuous
leafwise harmonic functions $\psi:M\rightarrow \mathbb{R}$    and $\Gamma$-equivariant continuous maps $\hat{\psi}:V\rightarrow \text{Har}(\mathbb{D})$. Furthermore, if $\Psi: M\rightarrow M_0$ is the harmonic
morphism induced from $\hat{\psi}$, then $\psi=\phi\circ\Psi$, and $\Psi$
is the unique morphism from $M$ to $M_0$ that satisfies this last equality.
\end{propo}

Therefore, if $V$ is any $\Gamma$-invariant  compact  subset of $X_0$
that contains non-constant functions, the foliated bundle with fiber $V$
will be an example of foliation for which the Liouville property fails, and
any example for which the property fails which is a foliated bundle over
$S$ is obtained in this way.

As an example, we note the  following alternative construction of the  foliation
of Theorem \ref{mainexample}.
   Consider  the element of $X_0$
whose boundary value is the function $f_0$ on $S^1$
such that 
\begin{equation*}
f_0(z)=
\begin{cases}
0 & \text{if }  \text{Im}(z)>0\\
1 & \text{if }  \text{Im}(z)\leq 0.
\end{cases}
\end{equation*}
Any   $g\in PSL(2,\mathbb{R})$ that fixes $f_0$ must also fix the
points $1, -1\in S^1$, which forces $g$ to be diagonal, as a simple
calculation shows. Therefore the orbit $X=PSL(2,\mathbb{R})\cdot f_0\in X_0$
 is identified with the quotient 
$PSL(2,\mathbb{R})/A$, where $A$ is the diagonal subgroup of $PSL(2,\mathbb{R})$.
(This quotient is the space of geodesics of the Poincar\'e disc.)
The closure of the orbit of $f_0$ is a topological sphere. 
Let $V$  be this orbit closure and $(M,\mathcal{F})$ the 
corresponding foliated bundle.  By noting that $f_0$ is the indicator function
of the arc $t\mapsto -e^{i\pi t}$, it is easy to check that the two descriptions of
the foliation of Theorem \ref{mainexample} are isomorphic. 

 \section{Relation  with  the  holomorphic case}
If the leaves of $(M,\mathcal{F})$ are K\"ahler manifolds, it makes sense to
ask whether a given continuous leafwise harmonic function
corresponds to the real part of a continuous leafwise holomorphic function.
This question is of interest since  the  holomorphically simple property seems
to be much easier to obtain  than its harmonic counterpart. In fact, we show
in \cite{ghani}  a number of results which prove
that a foliation is holomorphically simple, whose harmonic versions
are either not true, or
 we do not yet
know to hold. 

For example, if the leaves of $\mathcal{F}$
are complex manifolds (not necessarily K\"ahler) and the codimension is $1$,
then
 $(M,\mathcal{F})$ is holomorphically simple. It is still open whether a codimension
$1$ leafwise Riemannian  $(M,\mathcal{F})$ is harmonically simple.
It was also shown there that if the closure of each leaf contains
no more than countably many minimal sets, then $\mathcal{F}$ is holomorphically
simple. 
But the harmonic counterpart is not true, as the foliation of Theorem
\ref{mainexample} shows.

Regarding 
the relationship between 
holomorphic and harmonic functions, we limit ourselves here to making the following remark.  Recall that the (continuous)  tangential de Rham cohomology of a foliated
space is the cohomology of the complex $\Omega^*(M,\mathcal{F})$ of (continuous) tangentially smooth
differential forms, with the tangential exterior derivative. 
We denote the cohomology spaces by $H^*_{dR}(M,\mathcal{F})$.
If $f$ is a continuous leafwise harmonic function on $M$,
then we can obtain, on foliation boxes, continuous functions $\hat{f}$ which
are harmonic conjugates of $f$. These are well defined up to an additive constant,
so  the  various locally defined $d\hat{f}$ piece together to make a closed $1$-form  $\omega_f$ in $\Omega^1(M,\mathcal{F})$.  If the corresponding cohomology class, $[\omega_f]$,
is zero, then $f$ admits a global harmonic conjugate, $\hat{f}$,
and $f+i\hat{f}$ is a continuous leafwise holomorphic function. 
This remark yields the following proposition.

\begin{propo} Let $(M,\mathcal{F})$ be a foliation  of a compact
manifold $M$ by K\"ahler manifolds.
To each  leafwise harmonic continuous  function   $f:M\rightarrow \mathbb{R}$ 
is associated   a cohomology class  $[\omega_f]\in H^1_{dR}(M,\mathcal{F})$
having the property that $[\omega_f]=0$ if and only if $f$ is the real part of
a leafwise holomorphic continuous function.
\end{propo}

\begin{coro} 
Suppose that $H^1_{dR}(M,\mathcal{F})=\{0\}$,
for $(M,\mathcal{F})$ as in the previous proposition.  Then $(M,\mathcal{F})$ is harmonically
simple if and only if  it is holomorphically simple.
\end{coro}

We noted in Proposition \ref{holsimple} that the example given earlier of a
codimension $2$ foliated bundle
which is not harmonically simple is holomorphically simple. In particular,
it follows that it
has nontrivial first tangential cohomology. 
It is also interesting to observe that in codimension $1$ 
 only in very special cases does $H^1_{dR}(M,\mathcal{F})=\{0\}$ hold. 
For this reason,
 if it is at all  true in general
 that codimension $1$ foliations of a compact manifold by
Riemann surfaces are harmonically simple, we
expect the proof to be much harder 
than the related  holomorphic result   shown in \cite{ghani}.

\section{Dynamics  on $\text{Har}(\mathbb{D})$}
This section discusses some of the dynamical properties of the
action  of subgroups of $G=PSL(2,\mathbb{R})$ on $\text{Har}(\mathbb{D})$.
Unless further  hypothesis are explicitly 
assumed,  $\Gamma$ will denote an arbitrary subgroup of $G$. 
\subsection{General properties of $\Gamma$-actions on $\text{Har}(\mathbb{D})$}
We begin with a few general facts showing that limit sets of orbits in $\text{Har}(\mathbb{D})$
often contain constant functions. 
\begin{propo}
 Let  $g_n$  be  any unbounded sequence in    $PSL(2,\mathbb{R})$.
 Then for any $\varphi\in \text{Har}(\mathbb{D})$
whose boundary value is a continuous function
on $S^1$, all limit points  of  $\{g_n\cdot \varphi: n=1,2,\dots \}$ are constant functions.
\end{propo}
\begin{pf}
Write $g_n(z)=(\alpha_n z +\beta_n)/(\bar{\beta}_n z +\bar{\alpha}_n)$, where
$|\alpha_n|^2-|\beta_n|^2=1$, be
 an unbounded sequence of elements of $PSU(1,1)$. 
By taking a subsequence it can be assumed that   $\beta_n\neq 0$, and that
$\alpha_n/{\beta_n}$ and $\alpha_n/\bar{\beta_n}$ converge to points on the unit circle. Let $\eta$ be the limit of the latter sequence and observe that    
$ \varphi(g_n(z))=\varphi((\alpha_n/\bar{\beta}_n)(z+\beta_n/\alpha_n)/(
z+\bar{\alpha}_n/\bar{\beta}_n))\rightarrow \varphi(\eta)$, uniformly on $z$ in
any compact subset of $\mathbb{D}$.
\end{pf}

Since $PSL(2,\mathbb{R})$ has dimension $3$, we have the following.
\begin{coro}
Let $f\in C^0(S^1)$ and $\varphi$ the harmonic function on $\mathbb{D}$ with
boundary value $f$. Let $X$ denote the closure of the orbit 
$\Gamma\cdot \varphi$.
Then $X$ has topological dimension at most $3$  and
the complement of  $PSL(2,\mathbb{R})\cdot \varphi$ in $X$
consists of constant functions.
\end{coro}

As the previous proposition indicates, orbits in 
  $\text{Har}(\mathbb{D})$ 
with interesting dynamical properties are associated to harmonic functions
  whose boundary values   are not continuous functions on $S^1$.

Recall that $\text{Har}(\mathbb{D})^*$ denotes the nonconstant elements
of $\text{Har}(\mathbb{D})$.
\begin{propo}
$\Gamma$ fixes a point in 
  $\text{Har}(\mathbb{D})^*$ 
 if
and only if the $\Gamma$-action on $S^1$ is not ergodic.
\end{propo}
\begin{pf}
This is an immediate consequence of the
Poisson representation formula and
 the fact that the action of $\Gamma$ on $S^1$ is
not ergodic if and only if there exists a measurable, $\Gamma$-invariant, bounded   function on $S^1$.
\end{pf}
\begin{coro}
A lattice in $PSL(2,\mathbb{R})$ cannot have fixed points in $\text{Har}(\mathbb{D})^*$. 
\end{coro}

More generally,  we may ask whether a subgroup $\Gamma\subset PSL(2,\mathbb{R})$
leaves invariant a compact subset of $\text{Har}(\mathbb{D})^*$.
For cocompact lattices the answer is no.
\begin{propo}
A cocompact lattice in $PSL(2,\mathbb{R})$
does not leave invariant a compact subset of 
$\text{Har}(\mathbb{D})^*$.
\end{propo}
\begin{pf}
If such a $V$ existed, it would be possible to construct
a foliated space $(M,\mathcal{F})$ and a leafwise harmonic continuous functions
on $M$ which is not constant on any leaf. But this is clearly impossible.
\end{pf}

\subsection{Chaos}
Although there is no universally accepted mathematical definition of chaos,
one popular definition was proposed by Devaney in \cite{devaney}. He
isolates three properties as being the essential features of a chaotic dynamical
system: i) topological transitivity, ii) a dense set of periodic points, and iii)
sensitive dependence on initial conditions. The third condition
means the following. Let
$f:X\rightarrow X$ be a continuous map of a metric space $X$.
Then $f$ (or the $\mathbb{Z}$-action it generates on $X$) is said to 
have {\em sensitive dependence on initial conditions} if there exists
$\delta>0$ such that for every $x\in X$ and every neighborhood $N$ of $x$,
there exists $y\in N$ and a positive integer $n$ such that
$f^n(x)$ and $f^n(y)$ are  more than $\delta$ apart.

We show in this subsection   that   the
$\mathbb{Z}$-action generated by parabolic or hyperbolic elements
of $PSL(2,\mathbb{R})$ is chaotic. 
\begin{theor}
Let $\gamma$ be a hyperbolic or parabolic element of $PSL(2,\mathbb{R})$,
regarded as a transformation on $\text{Har}(\mathbb{D})$. Then $\gamma$ defines 
a chaotic dynamical system.
\end{theor}

It is proven in \cite{cairns} that  the first two conditions in the definition of chaos imply the 
third, so we only need to verify that these $\mathbb{Z}$-actions are topologically
transitive and have a dense set of periodic points.

Let $K_n=\{z\in \mathbb{C}: |z|\leq 1-1/n\}$
and define for $\phi\in \text{Har}(\mathbb{D})$,
\begin{equation*}
\|\phi\|=\sum_{n=1}^\infty  \frac{\sup_{z\in K_n}|\phi(z)|}{n^22^n}.
\end{equation*}
This norm induces a metric on 
$\text{Har}(\mathbb{D})$ compatible with the topology of uniform convergence on compact sets. Let  $f:S^1\rightarrow \overline{\mathbb{D}}$ be the boundary value
of $\phi$, which means that
\begin{equation*}
\phi(z)=\frac1{2\pi} \int_0^{2\pi}\frac{1-|z|^2}{|e^{i\theta}-z|^2}f(e^{i\theta})d\theta.
\end{equation*}
It is a simple calculation to check that
\begin{equation*}
\|\phi\|\leq \frac1{2\pi}\int_0^{2\pi}|f(e^{i\theta})|d\theta.
\end{equation*}

\begin{lem}\label{lem1}
The $\mathbb{Z}$-action   on $\text{Har}(\mathbb{D})$ generated by
any hyperbolic   element  $\gamma$ of  $  PSL(2,\mathbb{R})$ is topologically transitive.
\end{lem}
\begin{pf}
Without loss of generality, it can be assumed that $\gamma$ fixes $1$ and $-1$ in  $S^1$, and that $1$ is an expanding fixed point, while $-1$ is contracting.
Let $f_n$, $n=1,2,\dots$, be a weak*-dense sequence in the unit ball of $L^\infty(S^1)$.
It will be convenient to regard $S^1$ as the compactified real line, $\overline{\mathbb{R}}$, by means
of the map $\eta:z\mapsto i(1-z)/(1+z)$. Therefore, the $f_n$ are  now   regarded as functions
on the real line and $\gamma$ takes the form $\gamma(x)=\lambda x$, for 
some $\lambda>1$.

Choose a sequence of positive integers $k_1 < k_2<\dots$ having the property
\begin{equation*}
n(n+1)<\lambda^{k_{n+1}-k_n},
\end{equation*}
 for $n=1, 2,\dots$,  and define
$a_n=\lambda^{-k_n}/n$, $b_n=n\lambda^{-k_n}$. Notice that
the above inequality implies that $b_{n+1}<a_n$, so  
 the intervals $[a_n,b_n]$
are disjoint and accumulate at $0$.
Define now a function $f_{\infty}$ on $\overline{\mathbb{R}}$ by:
\begin{equation*}
f_{\infty}(x)=\begin{cases}f_n(\lambda^{k_n}x)&  \text{ for } x\in [-b_n, -a_n]\cup [a_n, b_n]\\
0 & \text{ for } x\in (-a_{n},-b_{n+1})\cup (b_{n+1}, a_n).
\end{cases}
\end{equation*}

Thus, we have by construction that    
$\gamma^{k_n}\cdot f_{\infty}:=f_\infty\circ \gamma^{-k_n}=f_n$ on  the set  
\begin{equation*}
A_n:= [-n, -1/n]\cup[1/n,n]=\gamma^{k_n}([-b_n, -a_n]\cup [a_n, b_n]).
\end{equation*}   
If $\varphi_{\infty}, \varphi_n$ are the elements of $\text{Har}(\mathbb{D})$ associated to
$f_\infty$, $f_n$ via the Poisson representation formula, then it
follows that  $\gamma^{k_n}\cdot \varphi_\infty$, $n=1, 2, \dots$ is 
arbitrarily close to $\varphi_n$  for large $n$. In fact,
let $B_n\subset S^1$ be the image of $A_n$ under $\eta^{-1}$. Note
that $f_n\circ \eta$ and $\gamma^{k_n}\cdot \varphi_\infty$ coincide on
$B_n$, whereas $B_n^c$ has length that goes to zero like $1/n$.
But $\|\gamma^{k_n}\cdot \varphi_\infty-\varphi_n\|\leq l(B_n^c)/\pi$
 so that $\mathbb{Z}\cdot \varphi_{\infty}$
is
dense
in $\text{Har}(\mathbb{D})$ as claimed. 
\end{pf}
\begin{lem} \label{lem2} The set of periodic points in
$\text{Har}(\mathbb{D})$ for the $\mathbb{Z}$-action generated by any 
     hyperbolic element $\gamma$ of $PSL(2,\mathbb{R})$
 is dense.
\end{lem}
\begin{pf}
For each $\varphi\in \text{Har}(\mathbb{D})$ and for each $\epsilon>0$,
we should find $\varphi_\epsilon\in \text{Har}(\mathbb{D})$ and
an integer $k\geq 2$ such that $\varphi_\epsilon\circ \gamma^k=\varphi_\epsilon$
and $\|\varphi-\varphi_\epsilon\| \leq \epsilon$. Let $f:S^1\rightarrow \overline{\mathbb{D}}$ be the boundary value of $\varphi$.

Let $B=B_n$ be as in the proof of the previous lemma, where $n$ is large enough 
so that  the length of the complement of $B$ satisfies $l(B^c)\leq \epsilon\pi$. Choose a positive integer $k$ such that $\gamma^k(B)$ and
$ B$ are disjoint and define 
\begin{equation*}
f_\epsilon=\sum_{m=-\infty}^{\infty}(\chi_{B}f)\circ\gamma^{mk}.
\end{equation*}
Notice that $f_\epsilon\circ \gamma^k=f_\epsilon$ and that $f_\epsilon$
still takes values into $\overline{\mathbb{D}}$. Furthermore,
\begin{align*}
\int_0^{2\pi}|f(e^{i\theta})-f_\epsilon(e^{i\theta}) |\ \! d\theta&=
\int_{B}|f(e^{i\theta})-f_\epsilon(e^{i\theta}) |\ \! d\theta
+ \int_{B^c}|f(e^{i\theta})-f_\epsilon(e^{i\theta}) |\ \! d\theta\\
&=\int_{B^c}|f(e^{i\theta})-f_\epsilon(e^{i\theta}) |\ \! d\theta\\
&\leq 2 l(B^c).
\end{align*}
Consequently, if $\varphi_\epsilon$ is the element of  $\text{Har}(\mathbb{D})$
with boundary value $f_\epsilon$, then $\varphi_\epsilon \circ \gamma^k=
\varphi_\epsilon$ and $\|\varphi-\varphi_\epsilon\|\leq \epsilon$.
  \end{pf}

If $\gamma$ is a parabolic element of $PSL(2,\mathbb{R})$, it
can be assumed (after conjugation) that its action on the extended real line
fixes $\infty$ so that $\gamma(x)=x+a$ for some   $a>0$.
The proofs of Lemmas \ref{lem1} and \ref{lem2} are easily modified
so as to apply to $\gamma$ parabolic. Notice that  now it
is convenient to take  $A_n=[-n,n]$.

\subsection{The conjugacy problem}
Given two distinct elements $\gamma_1$, $\gamma_2$ of 
  $PSL(2,\mathbb{R})$, it
is natural to ask whether or not they define dynamical systems 
 on $\text{Har}(\mathbb{D})$ that have the same topological dynamics.
More precisely,  can we find a homeomorphism $\Phi:\text{Har}(\mathbb{D})\rightarrow
\text{Har}(\mathbb{D})$ that conjugates (intertwines) the two transformations?

If $\gamma_1$ and $\gamma_2$ are both hyperbolic (resp., parabolic)
elements, then they are topologically conjugate. 
The conjugating homeomorphism  can be taken to be 
the map that associates to each   element of
$\text{Har}(\mathbb{D})$ with   boundary value $f$  
the element of $\text{Har}(\mathbb{D})$ with boundary value
$f\circ h$, where $h$ is a homeomorphism of $S^1$ that
conjugates $\gamma_1$ and $\gamma_2$, these two now regarded as transformations
of the circle. (It is easy to show that $h$ exists.)

It would be interesting to know, for example,  whether
a hyperbolic and a parabolic elements can be topologically conjugate. We
do not yet know the answer to this question.

There are many other 
topics about  the dynamics
of the action of $PSL(2,\mathbb{R})$ and its subgroups on $\text{Har}(\mathbb{D})$
which we did not consider here.  Some
of these  topics, particularly those about the ergodic theory of these actions,
will be taken up elsewhere.

\end{document}